\documentclass[a4, 10pt]{amsart}
\usepackage{amssymb}
\usepackage{amstext}
\usepackage{amsmath}
\usepackage{amscd}
\usepackage{latexsym}
\usepackage{amsfonts}

\theoremstyle{plain}
\newtheorem{thm}{Theorem}[section]
\newtheorem*{thm*}{Theorem}
\newtheorem*{cor*}{Corollary}

\newtheorem{prop}[thm]{Proposition}
\newtheorem{lem}[thm]{Lemma}

\newtheorem{cor}[thm]{Corollary}

\newtheorem*{claim*}{Claim}

\theoremstyle{definition}
\newtheorem{defn}[thm]{Definition}

\newtheorem{conj}[thm]{Conjecture}

\theoremstyle{remark}
\newtheorem*{pf}{{\sl Proof}}

\newtheorem*{tpf}{{\sl Proof of Theorem \ref{main}}}

\newtheorem*{ac}{\sc Acknowledgments}

\numberwithin{equation}{thm}
\def\Hom{\mathrm{Hom}}

\def\RHom{\mathrm{{\bf R}Hom}}
\def\Ext{\mathrm{Ext}}

\def\mod{\mathrm{mod}}

\def\m{\mathfrak m}

\def\p{\mathfrak p}

\def\D{\mathfrak D}

\def\depth{\mathrm{depth}}
\def\Supp{\mathrm{Supp}}

\def\Spec{\mathrm{Spec}}

\def\X{{\mathcal X}}

\def\G{{\mathcal G}}

\begin{document}

\title{On the number of indecomposable totally reflexive modules}
\author{Ryo Takahashi}
\address{Department of Mathematics, School of Science and Technology, Meiji University, 1-1-1 Higashimita, Tama-ku, Kawasaki 214-8571, Japan}
\email{takahasi@math.meiji.ac.jp}
\thanks{
\endgraf
}
\maketitle
\begin{abstract}
In this note, it is proved that over a commutative noetherian henselian non-Gorenstein local ring there are infinitely many isomorphism classes of indecomposable totally reflexive modules, if there is a nonfree cyclic totally reflexive module.
\end{abstract}
\section{Introduction}

Throughout this note, we assume that all rings are commutative and noetherian, and that all modules are finitely generated.

In the late 1960s, Auslander \cite{Auslander} introduced a homological invariant for modules which is called Gorenstein dimension, or G-dimension for short.
After that, he further developed the theory of G-dimension with Bridger \cite{AB}.
Many properties enjoyed by G-dimension are analogous to those of projective dimension.
An important feature is that G-dimension characterizes Gorenstein local rings exactly as projective dimension characterizes regular local rings.
A module of G-dimension zero is called a totally reflexive module.
Avramov and Martsinkovsky \cite{AM} and Holm \cite{Holm} proved that over a local ring any module $M$ of finite G-dimension admits an exact sequence $0 \to Y \to X \to M \to 0$ such that $X$ is totally reflexive and $Y$ is of finite projective dimension.
This result says that in the study of modules of finite G-dimension it is essential to consider totally reflexive modules.

On the other hand, Cohen-Macaulay local rings of finite Cohen-Macaulay type, namely Cohen-Macaulay local rings over which there are only finitely many isomorphism classes of indecomposable maximal Cohen-Macaulay modules, have been deeply studied since the work of Herzog \cite{Herzog}.
Under a few assumptions, Gorenstein local rings of finite Cohen-Macaulay type are hypersurfaces, and they have been classified completely.
Moreover, all isomorphism classes of indecomposable maximal Cohen-Macaulay modules over them are described concretely; see \cite{Yoshino} for the details.

Over a Gorenstein local ring, totally reflexive modules are the same as maximal Cohen-Macaulay modules.
Hence it is natural to expect that totally reflexive modules over an arbitrary local ring may behave similarly to maximal Cohen-Macaulay modules over a Gorenstein local ring, and we are interested in local rings over which there are only finitely many isomorphism classes of indecomposable totally reflexive modules; we want to determine all such isomorphism classes.
However, we guess that such ring cannot essentially exist in the non-Gorenstein case:

\begin{conj}\label{conj}
Let $R$ be a non-Gorenstein local ring.
Suppose that there is a nonfree totally reflexive $R$-module.
Then there are infinitely many isomorphism classes of indecomposable totally reflexive $R$-modules.
\end{conj}

In this conjecture, so as to exclude the case where all totally reflexive modules are free, it is assumed that there is a nonfree totally reflexive module.
Indeed, for instance, over a Cohen-Macaulay non-Gorenstein local ring with minimal multiplicity, every totally reflexive module is free \cite{Yoshino2}.

The author proved that the above conjecture is true over a henselian local ring of low depth:

\begin{thm}\cite{Takahashi0}\cite{Takahashi1}\cite{Takahashi2}\label{012}
Let $R$ be a henselian non-Gorenstein local ring of depth at most two.
Suppose that there is a nonfree totally reflexive $R$-module.
Then there are infinitely many isomorphism classes of indecomposable totally reflexive $R$-modules.
\end{thm}

The main purpose of this note is to prove that the conjecture is true over a henselian local ring having a nonfree cyclic totally reflexive module.

\begin{thm}\label{main}
Let $R$ be a henselian non-Gorenstein local ring.
Suppose that there is a nonfree cyclic totally reflexive $R$-module.
Then there are infinitely many isomorphism classes of indecomposable totally reflexive $R$-modules.
\end{thm}

This theorem says, for example, that if $R$ is a ring of the form $S[[X_1,\dots ,X_n]]/(f)$ where $S$ is a complete non-Gorenstein local ring and $f$ is a monomial, then there are infinitely many isomorphism classes of indecomposable totally reflexive $R$-modules.

In the next section, we will prove Theorem \ref{main} by using a theorem of Huneke and Leuschke \cite{HL} and Theorem \ref{012}.
In the last section, we will give several applications of Theorem \ref{main}.

\section{Proof of the theorem}

In this note, $(R, \m , k)$ is always a commutative noetherian local ring, and all $R$-modules are finitely generated.
We denote by $\mod\,R$ the category of finitely generated $R$-modules.
We begin with recalling the definition of a resolving subcategory.

\begin{defn}\label{resolv}
A full subcategory $\X$ of $\mod\,R$ is called {\it resolving} if the following hold.\\
(1) $\X$ contains $R$.\\
(2) $\X$ is closed under direct summands: if $M\in\X$ and $N$ is a direct summand of $M$, then $N\in\X$.\\
(3) $\X$ is closed under extensions: if there is an exact sequence $0 \to L \to M \to N \to 0$ in $\mod\,R$ with $L, N\in\X$, then $M\in\X$.\\
(4) $\X$ is closed under kernels of epimorphisms: if there is an exact sequence $0 \to L \to M \to N \to 0$ in $\mod\,R$ with $M, N\in\X$, then $L\in\X$.
\end{defn}

In this definition, the condition (3) especially says that $\X$ is closed under finite direct sums: if $M,N\in\X$, then $M\oplus N\in\X$.
Hence from (1) it follows that $\X$ contains all free $R$-modules.
Therefore, by (4), $\X$ is closed under syzygies: the (first) syzygy of any $R$-module in $\X$ is also in $\X$.

For an $R$-module $M$, we denote by $e(M)$ ($\nu (M)$, respectively) the multiplicity (the minimal number of generators, respectively) of $M$, namely,
$$
\begin{cases}
e(M)=\lim _{n\to\infty}\frac{d!}{n^d}\ell_R(M/\m^nM),\\
\nu (M)=\dim _k(M\otimes _Rk),
\end{cases}
$$
where $d=\dim M$ and $\ell _R(N)$ denotes the length of an $R$-module $N$.
Huneke and Leuschke essentially proved the following theorem in \cite[Theorems 1,3]{HL}.
(They actually proved the theorem in the case where $\X$ is the category of maximal Cohen-Macaulay $R$-modules.)

\begin{thm}[Huneke-Leuschke]
Let $\X$ be a full subcategory of $\mod\,R$ which is closed under extensions.
\begin{enumerate}
\item[{\rm (1)}]
Let $M, N\in\X$.
Assume that there are only finitely many isomorphism classes of $R$-modules $X\in\X$ with $e(X)=e(M)+e(N)$, and denote by $h$ the number of such isomorphism classes.
Then $\m ^h\Ext _R^1(M,N)=0$.
\item[{\rm (2)}]
Suppose that $\X$ is resolving.
Let $M\in\X$.
Assume that there are only finitely many isomorphism classes of indecomposable $R$-modules $X\in\X$ with $e(X)\leq\nu (M)\cdot e(R)$.
Then $M_\p$ is $R_\p$-free for any $\p\in\Spec\,R-\{\m\}$.
\end{enumerate}
\end{thm}

As a special case of the second assertion of this theorem, we obtain the following.

\begin{cor}\label{finitetype}
Let $\X$ be a resolving subcategory of $\mod\,R$.
Suppose that there are only finitely many isomorphism classes of indecomposable $R$-modules in $\X$.
Then $M_\p$ is $R_\p$-free for any $M\in\X$ and $\p\in\Spec\,R-\{\m\}$.
\end{cor}

Next, we recall the definition of a totally reflexive module.
Let $(-)^\ast$ be the $R$-dual functor $\Hom _R(-,R)$.

\begin{defn}
We say that an $R$-module $M$ is {\it totally reflexive} (or $M$ has {\it G-dimension zero}) if the natural homomorphism $M \to M^{\ast\ast}$ is an isomorphism and $\Ext _R ^i(M,R)=\Ext _R ^i(M^\ast ,R)=0$ for any $i>0$.
\end{defn}

We denote by $\G$ the full subcategory of $\mod\,R$ consisting of all totally reflexive $R$-modules.
Here, we state the properties of $\G$ which will be used later.

\begin{lem}\label{g}
\begin{enumerate}
\item[{\rm (1)}]
$\G$ is a resolving subcategory of $\mod\,R$.
\item[{\rm (2)}]
$\G$ is closed under $R$-dual, syzygies and finite direct sums.
\item[{\rm (3)}]
For any $M\in\G$, one has $\depth\,M=\depth\,R$.
\end{enumerate}
\end{lem}

\begin{pf}
(1) We refer to \cite[Lemma 2.3]{AM}, for example.

(2) It is easy to see from definition that if $M\in\G$ then $M^\ast\in\G$.
The remaining assertions follow from the arguments following Definition \ref{resolv}.

(3) See \cite[Proposition (4.12)]{AB} or \cite[Theorem (1.4.8)]{Christensen}.
\qed
\end{pf}

\begin{prop}\label{ess}
Suppose that there is a nonfree cyclic totally reflexive $R$-module $M$ such that $M_\p$ is $R_\p$-free for any $\p\in\Spec\,R-\{\m\}$.
Then $\depth\,R\leq 1$.
\end{prop}

\begin{pf}
Suppose that $\depth\,R\geq 2$.
We want to derive a contradiction.
We may assume $M=R/I$, where $I$ is an ideal of $R$ with $0\neq I\subseteq\m$.
Setting $J=(0:I)$, we have $J\neq R$, hence $I+J\subseteq\m$.
Dualizing the natural exact sequence $0 \to I \overset{\theta}{\to} R \to R/I \to 0$ and using that $R/I$ is assumed to be totally reflexive, gives an exact sequence
$$
0 \to J \to R \overset{\kappa}{\to} I^{\ast} \to \Ext ^1(R/I,R)=0,
$$
where $\kappa (1)=\theta$.
Thus we get an isomorphism $\lambda : R/J\to I^{\ast}$, where $\lambda (\overline{1})=\theta$.
Lemma \ref{g}(2) says that the $R$-modules $I$ and $I^\ast$ are totally reflexive.
Hence so is $R/J$, and there are isomorphisms
$$
I \to I^{\ast\ast} \overset{\lambda ^\ast}{\to} (R/J)^\ast \to (0:J).
$$
It is easy to check that the composite of these isomorphisms is an identity map; we obtain $I=(0:J)$.

Fix $\p\in\Spec\,R-\{\m\}$.
Since $(R/I)_\p$ is $R_\p$-free, one has either $IR_\p =0$ or $IR_\p =R_\p$.
If $IR_\p =0$, then $\p$ does not belong to $\Supp\,I$.
Noting that $\Supp\,I = V((0:I)) = V(J)$, we see that $J$ is not contained in $\p$.
If $IR_\p =R_\p$, then $I$ is not contained in $\p$.
This means that the ideal $I+J$ is $\m$-primary.
There is an exact sequence
$$
0 \to R/(I\cap J) \to R/I\oplus R/J \to R/(I+J) \to 0.
$$
Since the $R$-module $R/(I+J)$ has finite length, we have $\depth (R/(I+J))=0$.
According to Lemma \ref{g} parts (2) and (3), we get $\depth (R/I\oplus R/J)=\depth\,R\geq 2>0$.
Hence the depth lemma (cf. \cite[Proposition 4.3.1]{Roberts}) yields $\depth (R/(I\cap J))=1$.

Let $x\in I\cap J$.
Then, since $I=(0:J)$ and $J=(0:I)$, one has $xJ=xI=0$, which implies that $I+J\subseteq (0:x)$.
As $I+J$ is an $\m$-primary ideal, so is $(0:x)$.
Hence $\m ^r x=0$ for some $r>0$.
It follows that $I\cap J$ is an $R$-module of finite length.
Noting that $\depth\,R\geq 2>0$, one must have $I\cap J=0$.
Thus $2\leq\depth\,R=\depth (R/(I\cap J))=1$.
This contradiction proves the proposition.
\qed
\end{pf}

Now we can prove our main theorem.

\begin{tpf}
Suppose that $\G$ has only finitely many isomorphism classes of indecomposable $R$-modules.
Then Theorem \ref{012} implies that $\depth\,R\geq 3$.
But $M_\p$ is $R_\p$-free for any $M\in\G$ and $\p\in\Spec\,R-\{\m\}$ by Lemma \ref{g}(1) and Corollary \ref{finitetype}, hence $\depth\,R\leq 1$ by Proposition \ref{ess}.
This is a contradiction, which completes the proof of the theorem.
\qed
\end{tpf}

\section{Applications}

In this section, using Theorem \ref{main}, we give several results on the number of indecomposable totally reflexive modules.

\begin{cor}\label{0:x}
Let $(R, \m )$ be a henselian non-Gorenstein local ring.
If there exist $x,y\in\m$ such that $(0:x)=(y)$ and $(0:y)=(x)$, then there exist infinitely many nonisomorphic indecomposable totally reflexive $R$-modules.
\end{cor}

\begin{pf}
Noting that there is an exact sequence $\cdots\overset{x}{\to}R\overset{y}{\to}R\overset{x}{\to}\cdots$, we can easily check that the $R$-module $R/(x)$ is nonfree totally reflexive.
Hence the assertion follows from Theorem \ref{main}.
\qed
\end{pf}

\begin{cor}
Let $R$ be a complete non-Gorenstein local ring.
Then $S=R[[X_1,X_2,\dots ,X_n]]/(X_1^{a_1}X_2^{a_2}\cdots X_n^{a_n})$ admits infinitely many nonisomorphic indecomposable totally reflexive modules.
\end{cor}

\begin{pf}
Note that $S$ is faithfully flat over $R$.
Hence $S$ is also a complete non-Gorenstein local ring.
To show the corollary, we may assume $a_1>0$.
Then it is easily seen that $(0:_S X_1)=X_1^{a_1-1}X_2^{a_2}\cdots X_n^{a_n}S$ and $(0:_S X_1^{a_1-1}X_2^{a_2}\cdots X_n^{a_n})=X_1S$.
Thus we can apply Corollary \ref{0:x}.
\qed
\end{pf}

We denote by $\D (R)$ the derived category of $\mod\,R$.
Recall that an $R$-module $C$ is called {\it semidualizing} if the natural morphism $R\to\RHom _R(C,C)$ is an isomorphism in $\D (R)$, equivalently, the natural homomorphism $R\to\Hom _R(C,C)$ is an isomorphism and $\Ext _R ^i(C,C)=0$ for any $i>0$.
In the following, we consider the idealization $S=R\ltimes C$ of a semidualizing module $C$ over $R$.
There are two natural homomorphisms $\phi : R \to S$ and $\psi : S \to R$, which are given by $\phi (a)=(a,0)$ and $\psi (a,x)=a$.
Through the homomorphism $\phi$ ($\psi$, respectively), one can regard an $S$-module ($R$-module, respectively) as an $R$-module ($S$-module, respectively).
Note that through the composite of these homomorphisms the $R$-module structure is preserved since $\psi\phi$ is the identity map of $R$, but the $S$-module structure is not preserved in general.

\begin{lem}\cite[Lemma 3.2]{HJ}\label{hj}
Let $C$ be a semidualizing $R$-module, and set $S=R\ltimes C$.
Then there is a natural isomorphism $\RHom _R(-,C)\cong\RHom _S(-,S)$ of functors on $\D (R)$.
\end{lem}

Using this lemma, we can get the following result.
It says that a non-Gorenstein ring which is the idealization of a semidualizing module over a henselian local ring has infinitely many nonisomorphic totally reflexive modules.

\begin{cor}
Let $R$ be a henselian local ring, $C$ a semidualizing $R$-module, and $S=R\ltimes C$ the idealization.
Suppose that there are only finitely many nonisomorphic indecomposable totally reflexive $S$-modules.
Then $S$ is Gorenstein.
Hence $R$ is Cohen-Macaulay, and $C$ is a canonical module of $R$.
\end{cor}

\begin{pf}
The last assertion follows from \cite[Theorem 5.6]{FGR}.
(One can also prove it by using the isomorphism $\RHom _R (k,C)\cong\RHom _S(k,S)$ induced by Lemma \ref{hj}.)
Lemma \ref{hj} gives isomorphisms $\RHom _S(R,S)\cong C$ and $\RHom _S (C,S)\cong\RHom _R(C,C)\cong R$.
We easily see from these isomorphisms that $R$ is a totally reflexive $S$-module.
Note that $R$ is a nonfree cyclic $S$-module, and $S$ is henselian since $S$ is module-finite over $R$.
Therefore, Theorem \ref{main} implies that $S$ is Gorenstein.
\qed
\end{pf}

\begin{ac}
The author would like to thank Diana White for helpful comments.
\end{ac}



\begin{thebibliography}{99}


\bibitem{Auslander}
{\sc Auslander, M.}
Anneaux de Gorenstein, et torsion en alg\`{e}bre commutative.
S\'{e}minaire d'Alg\`{e}bre Commutative dirig\'{e} par Pierre Samuel, 1966/67. Texte r\'{e}dig\'{e}, d'apr\`{e}s des expos\'{e}s de Maurice Auslander, Marquerite Mangeney, Christian Peskine et Lucien Szpiro. \'{E}cole Normale Sup\'{e}rieure de Jeunes Filles
{\it Secr\'{e}tariat math\'{e}matique, Paris} 1967.

\bibitem{AB}
{\sc Auslander, M.}; {\sc Bridger, M.}
Stable module theory. 
Memoirs of the American Mathematical Society, No. 94
{\it American Mathematical Society, Providence, R.I.} 1969.

\bibitem{AM}
{\sc Avramov, L. L.}; {\sc Martsinkovsky, A.}
Absolute, relative, and Tate cohomology of modules of finite Gorenstein dimension.
{\it Proc. London Math. Soc. (3)} {\bf 85} (2002), no. 2, 393--440.

\bibitem{Christensen}
{\sc Christensen, L. W.}
Gorenstein dimensions.
Lecture Notes in Mathematics, 1747.
{\it Springer-Verlag, Berlin}, 2000. 

\bibitem{FGR}
{\sc Fossum, R. M.}; {\sc Griffith, P. A.}; {\sc Reiten, I.}
Trivial extensions of abelian categories.
Homological algebra of trivial extensions of abelian categories with applications to ring theory. Lecture Notes in Mathematics, Vol. 456.
{\it Springer-Verlag, Berlin-New York}, 1975.

\bibitem{Herzog}
{\sc Herzog, J.}
Ringe mit nur endlich vielen Isomorphieklassen von maximalen, unzerlegbaren Cohen-Macaulay-Moduln.
{\it Math. Ann.} {\bf 233} (1978), no. 1, 21--34.

\bibitem{Holm}
{\sc Holm, H.}
Gorenstein homological dimensions.
{\it J. Pure Appl. Algebra} {\bf 189} (2004), no. 1-3, 167--193.

\bibitem{HJ}
{\sc Holm, H}; {\sc J\o rgensen, P.}
Cohen-Macaulay homological dimensions.
preprint (2004), available from \texttt{http://www.maths.leeds.ac.uk/\~{}popjoerg/publications.html}.

\bibitem{HL}
{\sc Huneke, C.}; {\sc Leuschke, G. J.}
Two theorems about maximal Cohen-Macaulay modules.
{\it Math. Ann.} {\bf 324} (2002), no. 2, 391--404.

\bibitem{Roberts}
{\sc Roberts, P. C.}
Multiplicities and Chern classes in local algebra.
Cambridge Tracts in Mathematics, 133.
{\it Cambridge University Press, Cambridge}, 1998.

\bibitem{Takahashi0}
{\sc Takahashi, R.}
On the category of Gorenstein dimension zero.
{\it Math. Z.} {\bf 251} (2005), no. 2, 249--256.

\bibitem{Takahashi1}
{\sc Takahashi, R.}
On the category of modules of Gorenstein dimension zero. II.
{\it J. Algebra} {\bf 278} (2004), no. 1, 402--410.

\bibitem{Takahashi2}
{\sc Takahashi, R.}
Modules of G-dimension zero over local rings of depth two.
{\it Illinois J. Math.} {\bf 48} (2004), no. 3, 945--952.

\bibitem{Yoshino}
{\sc Yoshino, Y.}
Cohen-Macaulay modules over Cohen-Macaulay rings.
London Mathematical Society Lecture Note Series, 146.
{\it Cambridge University Press, Cambridge}, 1990. 

\bibitem{Yoshino2}
{\sc Yoshino, Y.}
Modules of G-dimension zero over local rings with the cube of maximal ideal being zero.
{\it Commutative algebra, singularities and computer algebra (Sinaia, 2002)}, 255--273, NATO Sci. Ser. II Math. Phys. Chem., 115, {\it Kluwer Acad. Publ., Dordrecht}, 2003. 


\end{thebibliography}
\end{document}